\documentclass[a4paper,12pt]{article}

\usepackage{amsfonts,amssymb,amsbsy,amsmath,amsthm,mathrsfs,enumerate,verbatim}
\usepackage{pdfsync}

\usepackage{graphics} 
\usepackage{epsfig} 
\usepackage{graphicx}
\usepackage{algorithmic,algorithm}
\usepackage{epstopdf}
\usepackage[colorlinks=true]{hyperref}
\hypersetup{urlcolor=blue, citecolor=red}

\usepackage{amsfonts}
\usepackage{amsmath}
\usepackage{amssymb}
\usepackage{graphicx}
\usepackage{epsf}

\usepackage{mathptmx}       
\usepackage{helvet}         
\usepackage{courier}        
\usepackage{type1cm}        
\usepackage{makeidx}         
\usepackage{graphicx}        
\usepackage{multicol}        
\usepackage[bottom]{footmisc}
\usepackage{amsfonts}
\usepackage{amsmath}
\usepackage{amssymb}
\usepackage{amsfonts}
\usepackage{graphicx}
\usepackage{epsf}
\usepackage{epsfig}

\usepackage{color}
\usepackage{afterpage}
\usepackage{verbatim}
\usepackage{algorithmic,algorithm}

\makeindex             

\graphicspath{
  {FIGURES/}
}

\begin{document}

\title{Numerical analysis of least squares and perceptron learning
  for classification problems}

\author{L.~Beilina  \thanks{Department of Mathematical Sciences, Chalmers University of Technology and University of Gothenburg, SE-42196 Gothenburg, Sweden, e-mail:
\texttt{\
  larisa@chalmers.se}}
}

%
%
\date{}
\maketitle 

\abstract{This work presents study on regularized and non-regularized
  versions of perceptron learning and least squares algorithms for
  classification problems.  The Fr\'{e}chet derivatives for
 least squares  and perceptron algorithms are
  derived. Different Tikhonov's regularization techniques for
  choosing the regularization parameter are discussed.  Numerical
  experiments demonstrate performance of perceptron and least squares
  algorithms to classify simulated and experimental data sets.  }

\section{Introduction}
\label{sec:intro}

 Machine learning is a field of artificial
intelligence which gives computer systems the ability to ``learn''
using available data. 
Recently machine learning algorithms become very popular for analyzing
of data and make prediction \cite{CB, hand, kots, MK}.
Linear models for classification is a part of supervised learning
\cite{CB}.  Supervised learning is machine learning task of learning a
function which transforms an input to an output data using available
input-output data.  In supervised learning, every example is a pair
consisting of an input object (typically a vector) and a desired
output value (also called the supervisory signal).  A supervised
learning algorithm analyzes the training data and produces an inferred
function, which can be used then for analyzing of new examples.
Supervised Machine Learning algorithms include linear and logistic
regression, multi-class classification, decision trees and support
vector machines.  In this work we will concentrate attention on study
 the linear and logistic regression algorithms.
Supervised learning problems are further divided into Regression and
Classification problems. Both problems have as goal the construction
of a model which can predict the value of the dependent
attribute from the attribute variables. The difference between these two
problems is the fact that the attribute is numerical for
regression and logical (belonging to class or not) for
classification.

 In this work are studied linear and polynomial classifiers, more
 precisely, the regularized versions of least squares and perceptron
 learning algorithms. The WINNOW algorithm for classification is also
 presented since it is used in numerical examples of Section
 \ref{sec:numex} for comparison of different classification
 strategies.  The classification problem is formulated as a
 regularized minimization problem for finding optimal weights in the
 model function. To formulate iterative gradient-based classification
 algorithms the Fr\'{e}chet derivatives for the non-regularized and
 regularized least squares algorithms are presented.  The Fr\'{e}chet
 derivative for the perceptron algorithm is also rigorously derived.

 Adding the regularization term in the functional leads to the optimal
 choice of weights such that they make a trade-off between observed
 data and obtaining a minimum of this functional.  Different rules are
 used for choosing the regularization parameter in machine learning,
 and most popular are early stopping algorithm, bagging and dropout
 techniques \cite{GBC}, genetic algorithms \cite{genetic}, particle
 swarm optimization \cite{pso1, pso2}, racing algorithms \cite{racing}
 and Bayesian optimization techniques \cite{bay1, bay2}. In this work
 are presented the most popular a priori and a posteriori Tikhonov's
 regularization rules for choosing the regularization parameter in the
 cost functional. 
Finally, performance of non-regularized versions of all classification
algorithms with respect to applicability, reliability and efficiency
is analyzed on simulated and experimental data sets \cite{iris,
  waves}.

The outline of the paper is as follows. In Section \ref{sec:2} are
briefly formulated non-regularized and regularized classification
problems.  Least squares for classification are discussed in Section
\ref{sec:3}. Machine learning linear and polynomial classifiers are
presented in Section \ref{sec:4}. Tikhonov's methods of regularization
for classification problems are discussed in Section
\ref{sec:5}. Finally, numerical tests are presented in Section
\ref{sec:numex}.

\section{Classification problem}

\label{sec:2}

The goal of regression is to predict the value of one or more
continuous target variables $t=\{t_i\}, i=1,...,m$ by knowing the values of input vector $x=\{x_i\}, i=1,...,m$. Usually, classification algorithms are working well for linearly separable  data sets.

\vspace{0.5cm}

\textbf{Definition}

Let $A$ and $B$ are two data sets of points in an $n$-dimensional
Euclidean space. Then $A$ and $B$ are linearly separable if there
exist $n + 1$ real numbers $\omega_1, ..., \omega_n, l$ such that
every point $x \in A$ satisfies $\sum _{i=1}^{n}\omega_{i} x_{i} > l$
and every point $x \in B$ satisfies $\sum _{i=1}^{n}\omega_{i} x_{i} <
-l$.

\vspace{0.3cm}

The classification problem is as follows:

\begin{itemize}
  
\item Suppose that we have data points   $\{x_i\}, i=1,...,m$ which
  are separated into two classes $A$ and $B$.  Assume that these
  classes are linearly separable.

\item Our goal is to find the decision line which will separate these
two classes. This line will also predict in which  class will
the new point fall.

\end{itemize}

In the non-regularized classification problem the goal is to find  optimal weights $\omega=(\omega_1,...,\omega_M)$, $M$ is the number  of weights,  in the functional
\begin{equation}\label{1}  
F(\omega) = \frac{1}{2}\| t - y(\omega) \|^2 = \frac{1}{2} \sum_{i=1}^m (t_i - y_i(\omega))^2
\end{equation}
with $m$ data points. Here, $ t =\{t_i\}, i = 1,...,m$, is the target function with known values, $ y(\omega) = \{ y_i(\omega)\} := \{ y(x_i, \omega)\}, i = 1,...,m$, is the classifiers model function.

{\small
  
  \begin{algorithm}[hbt!]
  \centering
  \caption{Gradient  Algorithm  for classification.\label{alg:ga}}
  \begin{algorithmic}[1]
    \STATE  Initialization:

    \begin{itemize}
 \item  Assume that every training example $\textbf{x}= (x_1,...,x_m)$  is described by $m$ attributes  with values $x_i=0$ or $x_i=1$.

\item Label examples of the first class with $t(\textbf{x})=1$  and  examples of the second class  $t(\textbf{x})=0$.

\item Denote by  $y(\textbf{x},\omega)$  the classifier's  model function.

\item Assume that all examples where $t(\textbf{x})=1$ are linearly
separable from examples where $t(\textbf{x})=0$.

\item Initialize weights $\omega^0 = \{\omega_i^0\}, i=1,...,M$ to small random numbers. Compute
  the sequence of ${\omega_i}^{N}$ for all $N >0$ in the following
  steps.
  
    \end{itemize}
    
    \STATE Compute  gradient
    
  \begin{equation}\label{5}
  G_i^{k} = -(t - y(\omega_i^k))\cdot  y'_{\omega_i}(\omega_i^k) + \gamma \omega_i^k, i=1,...,M.
  \end{equation}
  
  \STATE Update the unknown parameter $\omega := \omega^{k+1}$ using
  (\ref{5}) as
\begin{equation}\label{6}
  \omega_i^{k+1} =  \omega_i^{k} + \eta G^k,
\end{equation}
where $\eta$ is the learning rate or step size in the gradient update which is usually  taken as $\eta=0.5$ \cite{MK}.

  \STATE For the tolerance  $0 <\theta < 1$ chosen by the user, stop computing the functions
  $\omega_i^{k}$   and set $\omega_i^{N} = \omega_i^{k}$ if either
  $\|G(\omega^k)\|_{L_2} \leq \theta$, or  norms
   $\|G(\omega^k)\|_{L_2}$ abruptly grow, or computed
  $\|\omega^{k}\|_{L_2} $ are stabilized. Otherwise, set $k:=k+1$ and go to Step 2.
  \end{algorithmic}
\end{algorithm}
}

In the regularized classification problem to find optimal vector of weights $\omega = \{\omega_i\}, i=1,...,M$,
 to the functional \eqref{1} is added the regularization term  such that the functional is written as
  \begin{equation}\label{2}  
F(\omega) = \frac{1}{2}\| t - y(\omega) \|^2 + \frac{1}{2} \gamma \| \omega \|^2 = \frac{1}{2} \sum_{i=1}^m (t_i - y_i(\omega))^2  + \frac{1}{2} \gamma \sum_{j=1}^M |\omega_j|^2.
  \end{equation}
  Here, $\gamma$ is the regularization parameter, $\| \omega \|^2 = \omega^T \omega = \omega_1^2 + ... + \omega_M^2$, $M$ is the number  of weights.
  In order to find the optimal   weights in \eqref{1}  or in \eqref{2},
the following minimization problem should be solved
  \begin{equation}\label{minproblem}
\min_{\omega} F(\omega).
  \end{equation}
  Thus, we seek for a stationary point  of  \eqref{1} or \eqref{2} 
 with respect to $\omega$
such that
\begin{equation}\label{3}
F^{\prime }(\omega)(\bar{\omega})=0,
\end{equation}
where $F^{\prime }(\omega)$ is the Fr\'{e}chet derivative  acting on
$\bar{\omega}$.

 More precisely, for the functional \eqref{2}  
 we get
 \begin{equation}\label{4}
   \begin{split}
         F'(\omega)(\bar{\omega}) &=  \sum_{i=1}^M   F'_{\omega_i}(\omega)(\bar{\omega}_i), \\
  \frac{\partial F}{\partial \omega_i}(\omega)(\bar{\omega}_i) :=  F'_{\omega_i}(\omega)(\bar{\omega}_i) &=
   -(t - y)\cdot  y'_{\omega_i}(\bar{\omega}_i) + \gamma \omega_i(\bar{\omega}_i),
   ~i=1,..., M.
   \end{split}
 \end{equation}
 The  Fr\'{e}chet derivative of the functional \eqref{1} is obtained by taking $\gamma=0$ in \eqref{4}.
 To find optimal vector of weights $\omega = \{\omega_i\}, i=1,...,M$
 can be used the  \textbf{Algorithm \ref{alg:ga}} as well as   least squares or  machine learning algorithms.

For computation of the learning rate $\eta$ in the \textbf{Algorithm
  \ref{alg:ga}} usually is used optimal rule which can be derived
similarly as in \cite{person}. However, as a rule take $\eta=0.5$ in
machine learning classification algorithms \cite{MK}.  Among all other
regularization   methods applied in machine learning  \cite{bay1, racing, GBC, pso1, pso2, bay2, genetic}, the
regularization parameter $\gamma$ can be also computed using the Tikhonov's
theory for inverse and ill-posed problems by different algorithms
presented in \cite{BaK, BKS, BeK,IJ, KNS, TGSY}.  Some of these
algorithms are discussed in Section \ref{sec:5}.

\section{Least squares for classification}

\label{sec:3}

The linear regression is similar to the solution of linear least
squares problem  and can be used for classification problems appearing in
machine learning algorithms. We will revise solution of linear least squares
problem in terms of linear regression.

The simplest linear   model for regression is
\begin{equation}\label{linmod}
f(x,\omega) = \omega_0\cdot 1 + \omega_1 x_1 + ... + \omega_M x_M.
\end{equation}
Here, $\omega = \{ \omega_i \}, i=0,..., M$ are weights  with bias parameter $\omega_0$, $\{x_i\}, i=1,..., M$
are training examples.  Target values (known data) are $\{t_i\}, i=1,...,N$  which correspond to $\{x_i\}, i=1,...,M$.  Here, $M$ is the number of weights and
 $N$ is the number of data points.
The goal is to predict the value of $t$ in \eqref{1} for a
new value of $x$ in the model function \eqref{linmod}.

The linear model \eqref{linmod} can be written in the form
\begin{equation}\label{linmod2}
f(x,\omega) = \omega_0\cdot 1 + \sum_{i=1}^{M} \omega_i \varphi_i(x)  =  \omega_0 + \omega^T \varphi(x),
\end{equation}
where $\varphi_i(x), i=0,...,M$ are known basis functions with $\varphi_0(x) = 1$.

\subsection{Non-regularized least squares problem}

\label{sec:non-regls}

In non-regularized linear regression or least squares problem the goal
is  to minimize the sum of squares
\begin{equation}\label{linmod3}
  E(\omega) = \frac{1}{2}  \sum_{n=1}^{N}  (t_n - f(x,\omega))^2 = \frac{1}{2}  \sum_{n=1}^{N}  (t_n - \omega^T \varphi(x_n))^2 :=
  \frac{1}{2} \| t - \omega^T \varphi(x) \|_2^2
\end{equation}
to find
\begin{equation}\label{linmod4}
\min_{\omega} E(\omega) = \min_{\omega}  \frac{1}{2} \| t - \omega^T \varphi(x) \|_2^2.
  \end{equation}
The problem \eqref{linmod4} is a typical least squares problem of the
minimizing the squared residuals
\begin{equation}\label{linmod5}
\min_{\omega} \frac{1}{2} \| r(\omega) \|_2^2  = \min_{\omega}  \frac{1}{2} \| t - \omega^T \varphi(x) \|_2^2
  \end{equation}
with the residual $r(\omega) = t - \omega^T \varphi(x)$.
The test functions $\varphi(x)$ form the design matrix $A$
\begin{center}
\begin{equation}\label{linmod6}
    A=
\begin{bmatrix} 
1 & \varphi_1(x_1) & \varphi_2(x_1) & \hdots   & \varphi_M(x_1)\\ 
1 & \varphi_1(x_2) &  \varphi_2(x_2) &   \hdots &\varphi_M(x_2)\\
 1 & \varphi_1(x_3) & \varphi_2(x_3) & \hdots  &\varphi_M(x_3)\\
 \vdots & \vdots &  \ddots & \ddots &\vdots \\
 1 & \varphi_1(x_N) &  \varphi_2(x_N) &  \hdots & \varphi_M(x_N)
 \end{bmatrix},
\end{equation}
\end{center}
and the regression problem (or the least squares problem) is written as
\begin{equation}\label{linmod7}
 \min_{\omega} \frac{1}{2} \| r(\omega)\|_2^2  =  \min_{\omega} \frac{1}{2} \| A \omega - t \|_2^2,
\end{equation}
where $A$ is of the size $N \times M$ with $N > M$, $t$ is the target
vector  of the size $N$, and $\omega$ is vector of weights of the
size $M$.

 To
find minimum of the error function \eqref{linmod3} and derive the \emph{normal
  equations}, we look for the $\omega$ where the gradient of the norm
$\| r(\omega)\|_2^2 = ||A \omega - t||^2_2
= (A\omega - t)^T(A\omega - t)$ vanishes, or where $(\|r(\omega) \|_2^2)'_\omega =0$.
To derive the 
Fr\'{e}chet derivative,  we consider the difference
$\| r(\omega + e)\|_2^2 - \| r(\omega)\|_2^2$  and single out the linear part with respect to $\omega$. More precisely,  we get
\begin{equation*}
\begin{split}
  0  = &{\displaystyle\lim_{\|e\| \rightarrow 0}}\dfrac{(A(\omega+e)- t)^T(A(\omega + e) - t)-(A\omega- t)^T(A\omega - t)}{||e||_2} \\
  = &{\displaystyle\lim_{\|e\| \rightarrow 0}}
  \dfrac{((A \omega - t) + Ae )^T((A \omega  - t) + Ae)-(A\omega- t)^T(A\omega - t)}{||e||_2} \\
  =& {\displaystyle\lim_{\|e\| \rightarrow 0}} \dfrac{\|(A \omega - t) + Ae \|_2^2 - \| A\omega - t \|_2^2}{||e||_2} \\
  &=
{\displaystyle\lim_{\|e\| \rightarrow 0}}  \dfrac{\|A \omega - t\|_2^2   +
    2 (A \omega - t ) \cdot   Ae  + \| Ae \|_2^2 - \| A\omega - t \|_2^2}{||e||_2}  \\
   = & {\displaystyle\lim_{\|e\| \rightarrow 0}}\dfrac{2e^T(A^TA \omega - A^T t)+e^TA^TAe}{||e||_2}
\end{split}
\end{equation*}

Thus,
\begin{equation}\label{linmod8}
\begin{split}
 0 =  {\displaystyle\lim_{\|e\| \rightarrow 0}}\dfrac{2e^T(A^TA \omega - A^T t) + e^TA^TAe}{||e||_2}.
\end{split}
\end{equation}
The second term in \eqref{linmod8} can be estimated as
\begin{equation}
 {\displaystyle\lim_{\|e\| \rightarrow 0}} \frac{|e^TA^TAe|}{||e||_2} \leq {\displaystyle\lim_{\|e\| \rightarrow 0}} \frac{||A||_2^2 ||e||_2^2}{||e||_2}= {\displaystyle\lim_{\|e\| \rightarrow 0}} ||A||_2^2 ||e||_2 \to 0
\end{equation}
Thus, the first term in \eqref{linmod8}  must also be zero, or
\begin{equation}\label{linmod9}
A^T A \omega = A^T t
\end{equation}
Equations \eqref{linmod9}
is a  symmetric linear system of the $M \times M $
 linear equations for $M$ unknowns  called normal equations.

 Using  definition of the residual  in the functional
\begin{equation}
  \frac{1}{2} \|r(\omega) \|^2_2 = \frac{1}{2} \| A \omega - t \|_2^2
\end{equation}
 can  be computed the Hessian matrix $H = A^T A$.
If the Hessian matrix $H =  A^T A$ is positive definite, then $\omega$ is
indeed a minimum.

\textbf{Lemma}

The matrix  $A^T A$ is positive definite if
and only if the columns of $A$ are linearly independent, or when $rank(A)=M$ (full rank).

\emph{Proof}.

We have that $dim(A) = N \times M$, and thus, $dim(A^T A) = M \times M$.
Thus, $\forall v  \in R^M$ such that $v \neq 0$
\begin{equation}
v^T A^T A v = (A v)^T (Av) = \| Av\|_2^2 \geq 0. 
  \end{equation}
For positive definite matrix $A^TA$ we need to show that $v^T A^T A v > 0$.
Assume that $v^T A^T A v = 0$.
We observe that $A v = 0$ only if the linear combination  $\sum_{i=1}^M a_{ji} v_i = 0$. Here, $a_{ji}$ are elements of row $j$ in $A$.
This will be true only if columns of $A$ are linearly dependent or when $v=0$, but this is contradiction with assumption  $v^T A^T A v = 0$ since $v \neq 0$ and thus, the columns of $A$ are linearly independent and $v^T A^T A v > 0$. 
$\square$

The final conclusion is that if the matrix $A$ has a full rank
($rank(A)=M$) then the system (\ref{linmod9}) is of the size
$M$-by-$M$ and is symmetric positive definite system of normal
equations. It has the same solution $\omega$ as the least squares
problem $ \min_\omega \|A \omega - t\|^2_2$ and can be solved
efficiently via Cholesky decomposition \cite{BKK}.

\subsection{Regularized linear regression}

Let now the matrix $A$ will have entries $a_{ij} = \phi_j(x_i),
i=1,...,N; j = 1,...,M$.
Recall, that functions $\phi_j(x), j = 0,...,M$ are called basis functions which should be chosen and are known.
Then  the regularized least squares problem
takes the form
\begin{equation}\label{regleastsquares}
\min_{\omega} \frac{1}{2} \| r(\omega) \|_2^2  +  \frac{\gamma}{2}  \| \omega \|_2^2 = \min_{\omega}  \frac{1}{2} \|  A \omega - t \|_2^2 + \frac{\gamma}{2}  \| \omega \|_2^2.
\end{equation}

To minimize
the regularized squared residuals \eqref{regleastsquares}
we will again derive the \emph{normal
  equations}. Similarly as was derived the 
Fr\'{e}chet derivative for the  non-regularized regression problem \eqref{linmod7},    we look for the $\omega$ where the gradient of
$\frac{1}{2} ||A \omega - t||^2_2 + \frac{\gamma}{2}  \| \omega \|_2^2
= \frac{1}{2} (A\omega - t)^T(A\omega - t) + \frac{\gamma}{2}  \omega^T \omega$ vanishes. In other words, we consider
 the difference
$ (\| r(\omega + e)\|_2^2  + \frac{\gamma}{2}  \| \omega + e \|_2^2) - (\| r(\omega)\|_2^2 + \frac{\gamma}{2}  \| \omega \|_2^2)$, then   single out the linear part with respect to $\omega$ to obtain:
\begin{equation*}
\begin{split}
  0  = &\frac{1}{2}  {\displaystyle\lim_{\|e\| \rightarrow 0}}
  \dfrac{(A(\omega+e)- t)^T(A(\omega + e) - t)-(A\omega- t)^T(A\omega - t)}{||e||_2} \\
  &+ {\displaystyle \lim_{\|e\| \rightarrow 0}}  \dfrac{\frac{\gamma}{2} (\omega +e)^T(\omega + e) - \frac{\gamma}{2} \omega^T \omega}{||e||_2} \\
  =& \frac{1}{2} {\displaystyle\lim_{\|e\| \rightarrow 0}} \dfrac{\|(A \omega - t) + Ae \|_2^2 - \| A\omega - t \|_2^2}{||e||_2}  + {\displaystyle\lim_{\|e\| \rightarrow 0}}
 \dfrac{ \frac{\gamma}{2}(\| \omega + e\|_2^2  - \| \omega \|_2^2)}{ \|e\|_2}  \\
  &= \frac{1}{2}
{\displaystyle\lim_{\|e\| \rightarrow 0}}  \dfrac{\|A \omega - t\|_2^2   +
  2 (A \omega - t ) \cdot  Ae  + \| Ae \|_2^2 - \| A\omega - t \|_2^2}{||e||_2}   \\
+&
\frac{\gamma}{2} {\displaystyle\lim_{\|e\| \rightarrow 0}}  \dfrac{\|\omega \|_2^2   +
  2  e^T \omega   + \|e \|_2^2 - \| \omega \|_2^2}{||e||_2} \\
= & \frac{1}{2} {\displaystyle\lim_{\|e\| \rightarrow 0}}\dfrac{2e^T(A^TA \omega - A^T t)+e^TA^TAe}{||e||_2}
+ \frac{\gamma}{2}
 {\displaystyle\lim_{\|e\| \rightarrow 0}}\dfrac{2e^T \omega  + e^T e}{||e||_2}.
\end{split}
\end{equation*}

The  term 
\begin{equation}\label{linmodreg8}
\begin{split}
 {\displaystyle\lim_{\|e\| \rightarrow 0}} \frac{|e^TA^TAe|}{||e||_2} \leq  {\displaystyle\lim_{\|e\| \rightarrow 0}} \frac{||A||_2^2 ||e||_2^2}{||e||_2}= {\displaystyle\lim_{\|e\| \rightarrow 0}}  ||A||_2^2 ||e||_2 \to 0.
\end{split}
\end{equation}
Similarly, the  term 
\begin{equation}\label{linmodreg9}
\begin{split}
 {\displaystyle\lim_{\|e\| \rightarrow 0}} \frac{|e^T e|}{||e||_2} = {\displaystyle\lim_{\|e\| \rightarrow 0}} \frac{ ||e||_2^2}{||e||_2} \to 0.
\end{split}
\end{equation}
We finally get
 \begin{equation*}
\begin{split} 
  0 = {\displaystyle \lim_{\|e\| \rightarrow 0}}
  \dfrac{e^T(A^TA \omega - A^T t)}{||e||_2} + \dfrac{\gamma e^T \omega  }{||e||_2}.
\end{split}
 \end{equation*}
 The expression above  means that the factor
$A^TA \omega - A^T t  + \gamma \omega $  must also be zero, or
$$
(A^T A  + \gamma I) \omega = A^T t,
$$
where $I$ is the identity matrix.
This is a system of $M$ linear equations for $M$ unknowns, the normal equations
for regularized least squares.

\begin{figure}[hbt!]
 \begin{center}
   \begin{tabular}{cc}
{\includegraphics[scale=0.5, clip=true,]{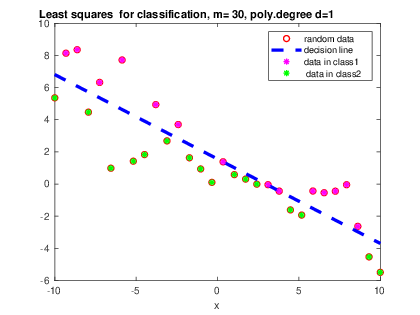}} &
  {\includegraphics[scale=0.5,  clip=true,]{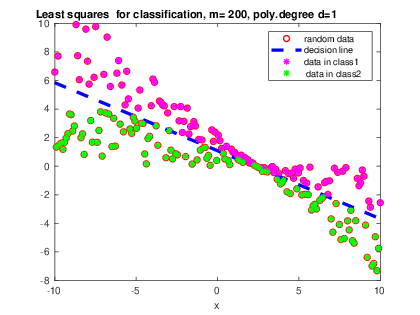}}  \\
\end{tabular}
 \end{center}
 \caption{Examples  of linear regression for classification  for different number of input points.}
  \label{fig:1}
\end{figure}

Figure  \ref{fig:1} shows that the linear regression or least squares minimization
$\min_\omega \| A \omega   -  t\|_2^2$ for classification is working fine when it is
known that two
classes are linearly separable. Here the linear model equation in the problem \eqref{linmod5} is
\begin{equation}\label{plmodel0}
 f(x,y,\omega) = \omega_0 + \omega_1 x + \omega_2 y
\end{equation}
and the target values of the vector $t=\{t_i\}, i=1,...,N$ in \eqref{linmod5} are
  \begin{equation}\label{plmodel1}
    t_i = \left \{
    \begin{array}{ll}
      1 & \textrm{red points}, \\
      0 &  \textrm{green points}. \\
    \end{array}
    \right.
    \end{equation}

The  elements of the design matrix  \eqref{linmod6} are given by
\begin{center}
\begin{equation}\label{plmodel2}
    A=
\begin{bmatrix} 
1 & x_1 & y_1 \\ 
1 & x_2 & y_2 \\
 1 & x_3 & y_3 \\
 \vdots & \vdots &  \ddots  \\
 1 & x_N &  y_N
 \end{bmatrix}.
\end{equation}
\end{center}

\subsection{Polynomial fitting to data in two-class model}

Let us consider the least squares classification in the two-class
model in the general case.  Let the first class consisting of $l$
points with coordinates $(x_i, y_i), i=1,...,l$ is described by it's
linear model
\begin{equation}\label{model1}
f_1(x, c) = c_{1,1} \phi_1(x) +  c_{2,1} \phi_2(x) + ... +  c_{n,1} \phi_n(x).
  \end{equation}
Let the second class  consisting of $k$ points  with coordinates  $(x_i, y_i), i=1,...,k$ is also described by the same linear model
\begin{equation}\label{model2}
f_2(x, c) = c_{1,2} \phi_1(x) +  c_{2,2} \phi_2(x) + ... +  c_{n,2} \phi_n(x).
  \end{equation}
Here, basis functions  are $\phi_j(x), j = 1,...,n$.  
 Our goal
 is to find the vector  of parameters $c= c_{i,1} = c_{i,2}, i=1,...,n$
 of the size $n$  which will fit best
to the data $y_i, i=1,...,m, m=k+l$    of both  model functions, $f_1(x_i, c), i=1,...,l$
and  $f_2(x_i, c), i=1,...,k$  with $f(x,c) = [f_1(x_i, c), f_2(x_i,c)]$
  such that the minimization problem
\begin{equation}\label{ch9_3}
\min_{c} \|  y - f(x, c) \|_2^2 = \min_{c} \sum_{i=1}^m ( y_i - f(x_i, c))^2
\end{equation}
is  solved with  $m=k+l$.
If the function $f(x, c)$  in \eqref{ch9_3}  is linear then we can reformulate
 the minimization problem \eqref{ch9_3}  as the following   least squares problem
\begin{equation}\label{ch9_4}
\min_{c} \| Ac -  y  \|_2^2, 
\end{equation}
where
the matrix $A$
 in the linear system
\begin{equation*}
  Ac=y
\end{equation*}
  will have entries $a_{ij} = \phi_j(x_i), i=1,...,m; j = 1,...,n$,
  i.e. elements of the matrix $A$ are created by basis functions
  $\phi_j(x), j = 1,...,n$.
  Solution of \eqref{ch9_4} can be found by the method of normal equations
  derived in  Section \ref{sec:non-regls}:
\begin{equation}\label{lsm}
 c = (A^T A)^{-1} A^T b = A^+ b
\end{equation}
with  pseudo-inverse matrix $A^+ := (A^T A)^{-1} A^T$.

\begin{figure}[hbt!]
 \begin{center}
   \begin{tabular}{cc}
{\includegraphics[scale=0.4, clip=true,]{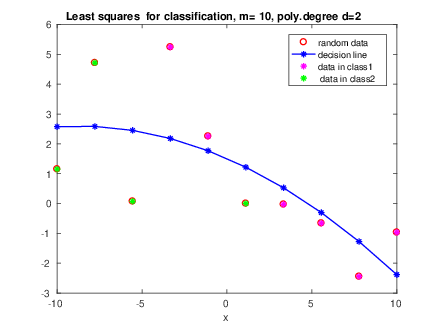}} &
{\includegraphics[scale=0.4,  clip=true,]{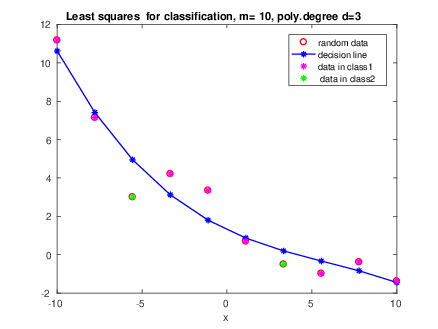}}  \\
{\includegraphics[scale=0.4, clip=true,]{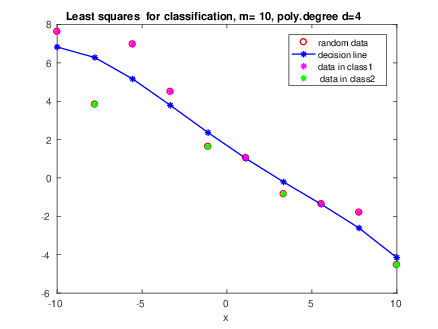}} &
{\includegraphics[scale=0.4,  clip=true,]{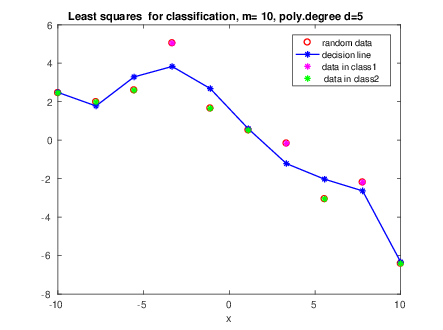}}  \\
\end{tabular}
 \end{center}
 \caption{Least squares in polynomial fitting to data for different degree of polynomial in the test functions \eqref{poly}.}
  \label{fig:3}
\end{figure}

For creating of elements of $A$ different basis functions can be chosen.
The polynomial test functions
\begin{equation}\label{poly}
\phi_j(x)   = x^{j-1},~~ j=1,...,n
\end{equation}
have been considered in the problem of fitting to a polynomial in examples presented in Figures \ref{fig:3}.
The
matrix $A$ constructed by these basis functions is a Vandermonde
matrix, and problems related to this matrix are
discussed in \cite{BKK}.
Linear splines (or hat functions)  and bellsplines
also can be used as basis functions \cite{BKK}.

Figures \ref{fig:3} present examples of polynomial fitting to data for
two-class model with $m=10$ using basis functions $\phi_j(x) =x^{j-1},
j=1,...,d$, where $d$ is degree of the polynomial.  Using these
figures we observe that least squares fit data well and even can
separate points in two different classes, although this is not always the case.
 Higher degree of polynomial separates two classes
better. However, since Vandermonde's matrix can be ill-conditioned for
high degrees of polynomial, we should carefully choose appropriate
polynomial to fit data.

{\small
\begin{algorithm}[hbt!]
  \centering
  \caption{Perceptron learning for classification.\label{alg:pl}}
  \begin{algorithmic}[1]
    \STATE  Initialization:

    \begin{itemize}
 \item  Assume that every training example $\textbf{x}= (x_1,...,x_n)$  is described by $n$ attributes.

\item Label examples of the first class with $c(\textbf{x})=1$  and  examples of the second class as $c(\textbf{x})=0$.

\item Let us denote by  $h(\textbf{x})$  the classifier's hypothesis  which  will have binary values $h(\textbf{x})=1$  or $h(\textbf{x})=0$.  Initialize $h(\textbf{x})=0$  for  examples  $c(\textbf{x})=1$  and  $h(\textbf{x})=1$  for  examples  $c(\textbf{x})=0$.
  

\item Assume that all examples of the first class where $c(\textbf{x})=1$ are linearly
separable from examples of the second class where $c(\textbf{x})=0$.

\item Initialize weights $\omega^0 = \{\omega_i^0\}, i=1,...,M$ to small random numbers. Compute
  the sequence of ${\omega_i}^{m}$ for all $m >0$ in the following
  steps.
  
    \end{itemize}
    
    \STATE  If $\sum_{i=0}^n \omega_i^m x_i > 0$  we will say that
    the example belongs to the first class and $h(x)=1$.

    \STATE If $\sum_{i=0}^n \omega_i^m x_i < 0$ we will say that the
    example belongs to the second class and $h(x)=0$.
    
  \STATE Update  weight $ \omega := \omega^{m+1} = \{\omega_i^{m+1} \}, ~i=1,...,M$ using
  
\begin{equation}\label{6}
\omega_i^{m+1} = \omega_i^m + \eta \cdot ([ c(\textbf{x})  - h(\textbf{x})]\cdot x_i + \gamma \cdot \omega_i^m),
\end{equation}
where $\eta$ is the learning rate usually  taken as $\eta=0.5$ \cite{MK}.

  \STATE  If $c(\textbf{x})=h(\textbf{x})$  for all learning examples - stop.
Otherwise   set $m := m+1$ and return to step 2.

  \end{algorithmic}
\end{algorithm}
}

\section{Machine learning linear and polynomial classifiers}

\label{sec:4}

\begin{figure}[hbt!]
 \begin{center}
   \begin{tabular}{cc}
{\includegraphics[scale=0.4, clip=true,]{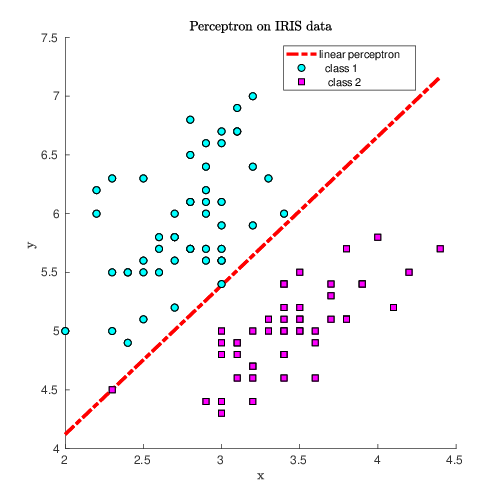}} &
  {\includegraphics[scale=0.4,  clip=true,]{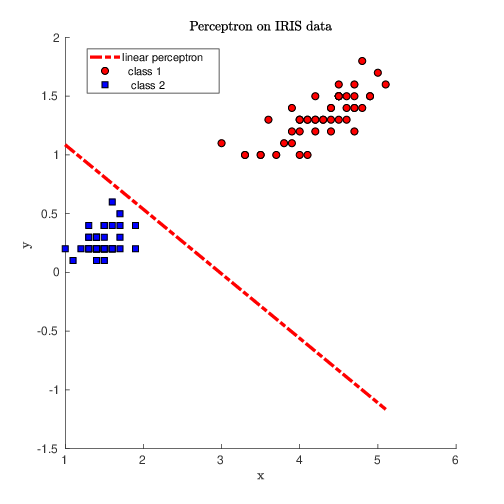}}  \\
\end{tabular}
 \end{center}
 \caption{Decision lines computed by the perceptron learning algorithm for separation of two  classes  using Iris dataset \cite{iris}.}
  \label{fig:4}
\end{figure}

In this section we will present the basic machine learning algorithms
for classification problems: perceptron learning and WINNOW
algorithms. Let us start with considering of an example: 
 determine the decision line for points presented in Figure
\ref{fig:4}.  One example on this figure is labeled as positive class, another one as
negative.  In this case, two classes are separated by the linear
equation with three weights $\omega_i, i=1,2,3$, given by
\begin{equation}\label{ml1}
\omega_1 + \omega_2 x + \omega_3 y = 0.
\end{equation}

In common case, two  classes can be separated by  the general equation
\begin{equation}
\omega_0 + \omega_1 x_1 + \omega_2 x_2 + ... + \omega_n x_n = 0
\end{equation}
which also can be written as
\begin{equation}
\omega^T x = \sum_{i=0}^n \omega_i x_i = 0
\end{equation}
with $x_0=1$.
If $n=2$ then the above equation defines a line,  if $n=3$ - plane, if $n>3$ - hyperplane.
The    problem is to determine weights $\omega_i$ and the task of machine
learning is to determine their appropriate values.  Weights $\omega_i,
i=1,...,n$ determine the angle of the hyperplane, $\omega_0$ is
called bias and determines the offset, or the hyperplanes distance from
the origin of the system of coordinates.

\subsection{Perceptron learning for classification}

The main idea of perceptron is binary classification.
The perceptron  computes a sum of weighted inputs
\begin{equation}\label{perc1}
y(x,\omega) = sign(\omega^T x)  = sign(\sum_{i=0}^n \omega_i x_i)
  \end{equation}
and uses then  binary classification.
When weights are computed, the linear classification boundary is defined by
\begin{equation*}
 \omega^T  x = 0.
  \end{equation*}
Thus, the perceptron algorithm determines weights $\omega_i, i=1,...,n$ in \eqref{perc1}
via binary classification.
 Binary classifier 
 decides whether or not an input $x$ belongs to some specific class:
 \begin{equation}\label{plmodel}
    sign(\omega^T x) = \left \{
    \begin{array}{ll}
      1, & \textrm{if $\sum_{i=1}^n \omega_i x_i + \omega_0 > 0$}, \\
      0, &  \textrm{otherwise},
    \end{array}
   \right.
    \end{equation}
 where $\omega_0$ is the bias. The bias does not depend on the input
 value $x$ and shifts the decision boundary. If the learning sets are
 not linearly separated the perceptron learning algorithm does not
 terminate and will never  converge and
 classify  data properly, see Figure \ref{fig:8}-a).
 
 The algorithm which determines weights in \eqref{perc1}  can be reasoned by minimization of the
   regularized residual
\begin{equation}\label{perc3}
  F(\omega) = \frac{1}{2}   \| r(x, \omega)\|_2^2 + \frac{1}{2} \gamma \| w \|_2^2 = \frac{1}{2} \| (t - y(x,\omega))\xi_\delta(x) \|_2^2
 + \frac{1}{2} \gamma \| w \|_2^2,
\end{equation}
 where $\xi_\delta(x)$ for a small $\delta$ is a data compatibility
 function to avoid discontinuities which can be defined similarly with
 \cite{IP} and $\gamma$ is the regularization parameter.  Taking $\gamma=0$
 algorithm will minimize the non-regularized residual \eqref{perc3}.
 Alternative,
 it can be minimized the residual
\begin{equation}\label{perc4}
  r(x, \omega) = - t^T y(x,\omega)  = - \sum_{i \in M} t_i y_i
\end{equation}
  over the set $M \subset
  \{ 1,..., m\} $ of the currently miss-classified patterns.

 The \textbf{ Algorithm  \ref{alg:pl}}  presents
the regularized perceptron learning algorithm where
in update of weights \eqref{6}    was used the following  regularized functional 
\begin{equation}\label{funcneuralp}  
  F(\omega) = \frac{1}{2}\| (t - y(x,\omega))\xi_\delta(x)  \|_2^2 + \frac{1}{2} \gamma \| w \|_2^2 = \frac{1}{2} \sum_{i=1}^m ((t_i - y_i(x,\omega))\xi_\delta(x) )^2 + \frac{1}{2} \gamma \sum_{i=1}^n  w_i^2.
 \end{equation}
Here, $\gamma$ is the regularization parameter,  $t$ is the target function, or class $c$ in the algorithm    \ref{alg:pl}, which takes values $0$ or $1$.

To find   optimal weights in \eqref{funcneuralp}  we need to solve the minimization problem  in the form \eqref{minproblem}
\begin{equation}
  F'(\omega)(\bar{\omega})=0,
  \end{equation}
where $F'(\omega)$ is a Frechet derivative acting on $\bar{\omega}$.
We will show how to derive  $F'(\omega)$ for
$y(x,\omega) = \omega^T x  = \sum_{i=0}^n \omega_i x_i$.
In a similar way can be derived  $F'(\omega)$   for $y(x,\omega)$  defined by
\eqref{perc1}.
We seek  the $\omega$ where the gradient of
$\frac{1}{2} ||r(x, \omega)||^2_2 + \frac{\gamma}{2}  \| \omega \|_2^2$ vanishes. In other words, we consider for \eqref{perc3}
 the difference
$ (\| r(x, \omega + e)\|_2^2  + \frac{\gamma}{2}  \| \omega + e \|_2^2) - (\| r(x, \omega)\|_2^2 + \frac{\gamma}{2}  \| \omega \|_2^2)$, then   single out the linear part with respect to $\omega$ to obtain:
\begin{equation*}
\begin{split}
  0  = &\frac{1}{2}  {\displaystyle\lim_{\|e\| \rightarrow 0}}
  \dfrac{ \|(t - y(x,\omega + e))\xi_\delta(x)\|_2^2 + \frac{\gamma}{2}  \| \omega + e \|_2^2  - \|(t - y(x,\omega))\xi_\delta(x)\|_2^2  -
   \frac{\gamma}{2}  \| \omega  \|_2^2 }{||e||_2} \\
 = &\frac{1}{2}  {\displaystyle\lim_{\|e\| \rightarrow 0}} \dfrac{ \|(t - \sum_{i=0}^n \omega_i x_i - \sum_{i=0}^n e_i x_i)\xi_\delta(x)\|_2^2  - \|(t - y(x,\omega))\xi_\delta(x)\|_2^2  }{||e||_2} \\ 
 &+ {\displaystyle \lim_{\|e\| \rightarrow 0}}  \dfrac{\frac{\gamma}{2} (\omega +e)^T(\omega + e) - \frac{\gamma}{2} \omega^T \omega}{||e||_2} \\
 = &\frac{1}{2}  {\displaystyle\lim_{\|e\| \rightarrow 0}} \dfrac{ \|(t -  y(x,\omega)  - e^T x)\xi_\delta(x)\|_2^2 - \|(t - y(x,\omega))\xi_\delta(x)\|_2^2}{||e||_2} \\ 
 &+ {\displaystyle \lim_{\|e\| \rightarrow 0}}  \dfrac{\frac{\gamma}{2} (\omega +e)^T(\omega + e) - \frac{\gamma}{2} \omega^T \omega}{||e||_2} \\
=& \frac{1}{2}  {\displaystyle \lim_{\|e\| \rightarrow 0}}
\dfrac{\|(t -  y(x,\omega))~\xi_\delta(x)\|_2^2 - 2(t - y(x,\omega))\cdot e^Tx~\xi_\delta(x) + \| e^T x~\xi_\delta(x)\|_2^2}{||e||_2}    \\
-& \frac{1}{2}  {\displaystyle \lim_{\|e\| \rightarrow 0}}\dfrac{\|(t -  y(x,\omega))~\xi_\delta(x)\|_2^2}{||e||_2} 
  + {\displaystyle \lim_{\|e\| \rightarrow 0}}
  \dfrac{\frac{\gamma}{2} (\omega +e)^T(\omega + e) - \frac{\gamma}{2} \omega^T \omega}{||e||_2} \\
 =& \frac{1}{2}  {\displaystyle\lim_{\|e\| \rightarrow 0}} \dfrac{ - 2 (t - y(x,\omega))\cdot e^Tx~ \xi_\delta(x) + \| e^T x~\xi_\delta(x)\|_2^2}{||e||_2}   + \frac{\gamma}{2}
      {\displaystyle\lim_{\|e\| \rightarrow 0}}\dfrac{2e^T \omega  + e^T e}{||e||_2}.
\end{split}
\end{equation*}
The second part in the last  term of the above expression  is estimated as in \eqref{linmodreg9}.
The second part in the first term is estimated as  
\begin{equation}\label{linmodreg8}
\begin{split}
  &{\displaystyle\lim_{\|e\| \rightarrow 0}} \frac{|(e^T x~\xi_\delta(x))^T~e^T x~\xi_\delta(x) |}{||e||_2} =
 {\displaystyle\lim_{\|e\| \rightarrow 0}} \frac{|(x^T e ~\xi_\delta(x))^T~ x^T e ~\xi_\delta(x) |}{||e||_2} \\ 
  &\leq  {\displaystyle\lim_{\|e\| \rightarrow 0}} \frac{|| x ~\xi_\delta(x)||_2^2 ||e||_2^2}{||e||_2}= {\displaystyle\lim_{\|e\| \rightarrow 0}}  ||x ~\xi_\delta(x) ||_2^2 ||e||_2 \to 0.
\end{split}
\end{equation}

We finally get
 \begin{equation*}
\begin{split} 
  0 = {\displaystyle \lim_{\|e\| \rightarrow 0}}
  -\dfrac{x^T e( t - y(x,\omega)) ~\xi_\delta(x)}{||e||_2} + \dfrac{\gamma e^T \omega  }{||e||_2}.
\end{split}
 \end{equation*}
 The expression above  means that the factor
$- x^T ( t - y(x,\omega)) ~\xi_\delta(x)  + \gamma \omega $  must also be zero, or
   \begin{equation}\label{frechet2p}
     \begin{split}
   F'(\omega)(\bar{\omega}) &=  \sum_{i=1}^n   F'_{\omega_i}(\omega)(\bar{\omega}_i), \\
   F'_{\omega_i}(\omega)(\bar{\omega}_i) &=
-(t - y) \cdot \xi_\delta(x)\cdot  y'_{\omega_i}(\bar{\omega}_i) + \gamma \omega_i = -(t - y) \cdot x_i \cdot \xi_\delta(x_i) + \gamma \omega_i, ~i =1,...,n.
   \end{split}
 \end{equation}
   The non-regularized version of the perceptron
   \textbf{Algorithm \ref{alg:pl}} is obtained  taking  $\gamma=0$ in \eqref{frechet2p}.

\subsection{Polynomial of the second order}

Coefficients of polynomials of the second order can be obtained by the same technique as coefficients for linear classifiers.
The second order polynomial function is:
\begin{equation}\label{secondorder}
  \omega_0 + \omega_1 \underbrace{x_1}_{z_1} + \omega_2 \underbrace{x_2}_{z_2}
  +  \omega_3 \underbrace{x_1^2}_{z_3} + \omega_4 \underbrace{x_1 x_2}_{z_4}
  + \omega_5 \underbrace{x_2^2}_{z_5} = 0.
\end{equation}

This polynomial can be converted to the linear classifier if we introduce notations:
$$
z_1 = x_1, z_2=x_2, z_3= x_1^2, z_4 =  x_1 x_2, z_5= x_2^2.
$$

Then equation \eqref{secondorder}  can be written  in new variables as
\begin{equation}\label{secondorder2}
  \omega_0 + \omega_1 z_1 + \omega_2 z_2 +  \omega_3 z_3 + \omega_4 z_4 +
  \omega_5 z_5 = 0
\end{equation}
which is already the linear   function.  Thus, the Perceptron learning
 \textbf{ Algorithm  \ref{alg:pl}} can be used to determine weights  $\omega_0, ...,  \omega_5$ in
\eqref{secondorder2}.

 Suppose that  weights  $\omega_0, ...,  \omega_5$ in
 \eqref{secondorder2}  are computed. To present the decision  line one
 need to solve the  following quadratic equation  for  $x_2$: 
\begin{equation}\label{secondorder3}
\omega_0 + \omega_1 x_1 + \omega_2 x_2 +  \omega_3 x_1^2 + \omega_4 x_1 x_2 + \omega_5 x_2^2 = 0
\end{equation}
with known weights  $\omega_0, ...,  \omega_5$  and known  $x_1$
which can be rewritten as
\begin{equation}\label{secondorder4}
\underbrace{\omega_5}_a x_2^2    +  x_2 \underbrace{( \omega_2  + \omega_4 x_1)}_b  + \underbrace{\omega_0 + \omega_1 x_1 + \omega_3 x_1^2}_c =0,
\end{equation}
or in the form
\begin{equation}\label{secondorder5}
  a x_2^2    + b x_2   +  c =0
\end{equation}
with known coefficients  $ a= \omega_5, b = \omega_2 + \omega_4 x_1, c = \omega_0 + \omega_1 x_1 + \omega_3 x_1^2$.
Solutions of \eqref{secondorder5} will be
\begin{equation}\label{secondorder6}
  \begin{split}
x_2 &= \frac{-b \pm \sqrt{D}}{2a}, \\
D &= b^2 -4ac.
\end{split}
\end{equation}

Thus, to present the decision line for polynomial of the second order,
first should be computed weights $\omega_0, ..., \omega_5$, and then
the quadratic equation \eqref{secondorder4}   should be solved the 
solutions of which are
given by \eqref{secondorder6}. Depending on the classification problem
and set of admissible parameters for classes, one can then decide which
one classification line should be presented, see examples in section
\ref{sec:numex}.

\begin{figure}[hbt!]
 \begin{center}
   \begin{tabular}{cc}
{\includegraphics[scale=0.5, clip=true,]{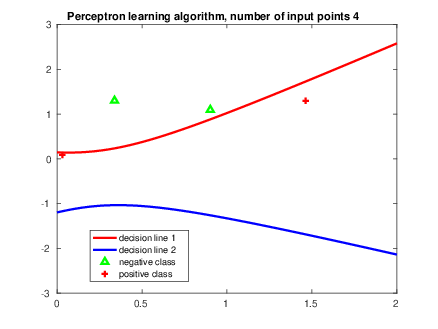}} &
{\includegraphics[scale=0.5, clip=true,]{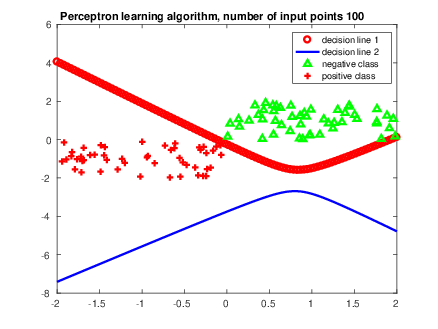}} 
\end{tabular}
 \end{center}
 \caption{Perceptron learning algorithm for separation of two classes  by polynomials of the second order.}
  \label{fig:5}
\end{figure}

\begin{algorithm}[hbt!]
  \centering
  \caption{WINNOW for classification \label{alg:wc}}
  \begin{algorithmic}[1]

    \STATE   Initialization:


    \begin{itemize}
 \item  Assume that every training example $\textbf{x}= (x_1,...,x_n)$  is described by $n$ attributes.

\item Label examples of the first class with $c(\textbf{x})=1$  and  examples of the second class as $c(\textbf{x})=0$. Assume that all examples of the first class where $c(\textbf{x})=1$ are linearly
separable from examples of the second class where $c(\textbf{x})=0$.

\item  Initialize  the classifier's hypothesis $h(\textbf{x})=0$  for  examples  $c(\textbf{x})=1$  and  $h(\textbf{x})=1$  for  examples  $c(\textbf{x})=0$.
  

\item Choose parameter $\alpha >1$, usually $\alpha=2$.
  
\item Initialize weights $\omega^0 = \{\omega_i^0\}, i=1,...,M$ to small random numbers. Compute
  the sequence of ${\omega_i}^{m}$ for all $m >0$ in the following
  steps.
  
    \end{itemize}

\STATE    If $\sum_{i=0}^n \omega_i^m x_i > 0$  we will say that the example is positive and $h(x)=1$.

\STATE If $\sum_{i=0}^n \omega_i^m x_i < 0$ we will say the the example is negative and $h(x)=0$.

\STATE  Update every weight using the formula
$$
\omega_i^{m+1} = \omega_i^m \cdot \alpha^{ (c(\textbf{x})  - h(\textbf{x}))\cdot x_i }.
$$

\STATE  If $c(\textbf{x})=h(\textbf{x})$  for all learning examples - stop.
Otherwise set $m:=m+1$ return to step 1.

  \end{algorithmic}
\end{algorithm}

\subsection{WINNOW learning algorithm}

To be able compare perceptron with
other machine learning algorithms,  we present here one more learning
algorithm which is very close to the perceptron and called WINNOW.  Here is described the simplest version of this
algorithm without regularization.  The regularized version of WINNOW
is analyzed in \cite{zhang}.  Perceptron learning algorithm uses
additive rule in the updating weights, while WINNOW algorithm  uses
multiplicative rule: weights are multiplied in this rule.  The WINNOW
algorithm \textbf{ Algorithm  \ref{alg:wc}} is written for $c=t$ and $y=h$ in \eqref{frechet2p}.  We
will again assume that all examples where $c(\textbf{x})=1$ are
linearly separable from examples where $c(\textbf{x})=0$.

\section{Methods of Tikhonov's  regularization    for classification problems}

\label{sec:5}

To solve the regularized classification problem the regularization
parameter $\gamma$ can be chosen by the same methods which are used
for the solution of ill-posed problems. For different Tikhonov's
regularization strategies we refer to \cite{BaK, BKS, BeK, IJ, KNS,
  TGSY, TA}.  In this section we will present main
methods of Tikhonov's regularization
 which follows ideas of
\cite{BaK, BKS, BeK, IJ,TGSY,TA}.

\textbf{Definition} Let $B_{1}$ and $B_{2}$ be two Banach
spaces and $G\subset B_{1}$ be a set. Let  $y:G\rightarrow B_{2}$
be one-to-one. Consider the equation
\begin{equation}
y( \omega) = t, \label{1.28}
\end{equation}
where $t$ is the target function and $y(\omega)$  is the model function in the classification problem.
Let $t^{\ast }$ be the  noiseless target function in equation (\ref%
{1.28}) and $\omega^{\ast }$ be the ideal noiseless weights corresponding to
$
t^{\ast }, ~y(\omega^{\ast }) =  t^{\ast}$.  For every $\delta \in \left(
0,\delta _{0}\right), ~\delta _{0}\in \left(
0,1\right)  $ denote
\begin{equation*}
K_{\delta }( t^{\ast }) = \left\{ z\in B_{2}:\left\Vert z- t^{\ast
}\right\Vert _{B_{2}}\leq \delta \right\}.
\end{equation*}

Let $\gamma >0$ be a parameter and $R_{\gamma}:K_{\delta _{0}}(t^{\ast}) \rightarrow G$ be a continuous operator depending on the
parameter $\gamma$. The operator $R_{\gamma }$ is called the  \emph{%
    regularization operator}  for \eqref{1.28}
if there exists a
function $\gamma_0 \left( \delta \right) $ defined for $\delta \in \left(
0,\delta _{0}\right) $ such that
\begin{equation*}
\lim_{\delta \rightarrow 0} \left\Vert R_{\gamma_0 \left( \delta \right)
}\left( t_{\delta }\right) - \omega^{\ast }\right\Vert _{B_{1}}=0.
\end{equation*}
The parameter $\gamma$ is called the regularization parameter and the
procedure of constructing the approximate solution $\omega_{\gamma(
  \delta)} = R_{\gamma(\delta) }(t_{\delta })$ is called the
\emph{regularization procedure}, or simply \emph{regularization}.

One can use different regularization procedures for the same
classification problem. The regularization parameter $\gamma$  can
be even the vector of regularization parameters depending on number of
iterations in the used classification method, the tolerance chosen by the user, number of classes, etc..

For two Banach spaces $B_{1}$ and $B_{2}$  let $Q$ be another
space, $Q\subset B_{1}$ as a set and $\overline{Q}=B_{1}$. In
addition, we assume that $Q$ is compactly embedded in $B_{1}.$ Let
$G\subset B_{1}$\ be the closure of an open set$.$ Consider a
continuous one-to-one function $y:G\rightarrow B_{2}.$ Our goal is to
solve
\begin{equation}
y( \omega) = t,~~ \omega\in G.  \label{1.37}
\end{equation}%
Let
\begin{equation}
y(\omega^*) = t^*,\quad \left\Vert t - t^{\ast }\right\Vert
_{B_{2}}<\delta .  \label{1.38}
\end{equation}%
To find an approximate solution of equation (\ref{1.37}), we  construct the
Tikhonov regularization functional  $J_{\gamma}(\omega),$
\begin{equation}
J_{\gamma}(\omega) =\frac{1}{2}\left\Vert y(\omega) - t \right\Vert
_{B_{2}}^{2}+ \frac{\gamma}{2} \| \omega \|_{B_1}^2 := \varphi(\omega) +
\frac{\gamma}{2} \psi(\omega),
\label{1.39}
\end{equation}%
\begin{equation*}
J_{\gamma}: G\rightarrow \mathbb{R},
\end{equation*}
where $\gamma = \gamma(\delta) >0$ is a small regularization
parameter.

The regularization term $\frac{\gamma}{2} \psi(\omega) $ encodes
a priori available information about the unknown solution such that sparsity,
smoothness, monotonicity, etc...  Regularization  term $ \frac{\gamma}{2} \psi(\omega) $  can be chosen
in different norms, for example:

  \begin{itemize}

  \item  $ \frac{\gamma}{2} \psi(\omega) = \frac{\gamma}{2} \|  \omega \|^p_{L^p}, ~~ 1 \leq  p \leq 2$.

  \item  $\frac{\gamma}{2} \psi(\omega) = \frac{\gamma}{2} \| \omega \|_{TV}$,  TV means ``total variation''.

  \item $ \frac{\gamma}{2} \psi(\omega) = \frac{\gamma }{2} \| \omega \|_{BV}$, BV   means ``bounded variation'', a real-valued function whose TV is bounded (finite).

    \item  $ \frac{\gamma}{2} \psi(\omega) = \frac{\gamma }{2} \|\omega\|_{H^1}^2$.

      \item $ \frac{\gamma}{2} \psi(\omega) =  \frac{\gamma }{2}( \|\omega \|_{L^1}  + \| \omega \|^2_{L^2})$.

  \end{itemize}

We   consider the following Tikhonov functional for
regularized classification problem 
\begin{equation}
J_{\gamma}(\omega) =\frac{1}{2}\left\Vert y(\omega) - t \right\Vert
_{L_2}^{2}+ \frac{\gamma}{2} \| \omega - \omega_0\|_{L_2}^2 := \varphi(\omega) +
\frac{\gamma}{2} \psi(\omega),
\label{1.39}
\end{equation}
where terms $\varphi(\omega), \psi(\omega)$ are considered  in $L_2$ norm which is the classical Banach space.
In \eqref{1.39}  $\omega_{0}$ is a good first
approximation for the exact   weight function $\omega^{\ast }$, which is
 called also the first guess or the  first approximation. For discussion about how the first guess in the functional \eqref{1.39}  should be chosen we refer to \cite{BaK, BeK, KBB}.

In this section we will discuss  following rules for choosing regularization parameter in \eqref{1.39}:
  
\begin{itemize}
  
\item A-priori rule (Tikhonov's regularization)
  
  \begin{itemize}

  \item   For $\| t - t^*\| \leq \delta$ a priori rule requires (see details in \cite{BaK, BKS}):
\begin{equation*}
  \lim_{\delta \rightarrow  0} \gamma(\delta) \to 0,  ~ \lim_{\delta \rightarrow 0} \frac{\delta^2}{\gamma(\delta)} \to 0.
\end{equation*}

  \end{itemize}

 \item A-posteriori rules:
  
  \begin{itemize}

  \item  Morozov's discrepancy principle  \cite{IJ, M, TGSY}.

  \item Balancing principle \cite{IJ}.

    \end{itemize}

\end{itemize}

A-priori rule and Morozov's discrepancy are most popular methods for the case when there
exists estimate of the noise level $\delta$ in data $t$. Otherwise it is
  recommended to use   balancing principle or other a-posteriori rules presented in \cite{BaK, BKS, BeK, IJ, KNS,
  TGSY, TA}.

\subsection{The Tikhonov's regularization }

The goal of regularization is 
to construct sequences $\left\{ \gamma\left( \delta _{k}\right) \right\}
,\left\{ \omega_{\gamma \left( \delta _{k}\right) }\right\} $ in a stable way so that
\begin{equation*}
\lim_{k\rightarrow \infty }\left\Vert \omega_{\gamma \left( \delta _{k}\right)
}- \omega^{\ast }\right\Vert _{B_{1}}=0,
\end{equation*}
where a sequence $\left\{ \delta _{k}\right\} _{k=1}^{\infty }$
is such that
\begin{equation}\label{convseq}
\delta _{k}>0,~~\lim_{k\rightarrow \infty }\delta _{k}=0.
\end{equation}

Using (\ref{1.38}) and (\ref{1.39}), we obtain
\begin{align}
  J_{\gamma}\left( \omega^{\ast }\right)
  &= \frac{1}{2}\left\Vert y(\omega^{\ast }) -  t \right\Vert
_{B_{2}}^{2}+\frac{\gamma(\delta) }{2}\left\Vert \omega^{\ast } - \omega_{0}\right\Vert _{Q}^{2} \\
&\leq
\frac{\delta ^{2}}{2} + \frac{\gamma(\delta)}{2}\left\Vert \omega^{\ast } - \omega_{0}\right\Vert _{Q}^{2}.
\label{1.40}
\end{align}
Let
\begin{equation*}
m_{\gamma \left( \delta _{k}\right) }=\inf_{G}J_{\gamma \left( \delta
_{k}\right) }\left( \omega\right) .
\end{equation*}%
By (\ref{1.40})
\begin{equation*}
  m_{\gamma \left( \delta _{k}\right) } \leq \frac{\delta _{k}^{2}}{2}+\frac{%
\gamma \left( \delta _{k}\right) }{2}\left\Vert \omega^{\ast } - \omega_{0}\right\Vert_{Q}^{2}.
\end{equation*}
Hence, there exists a point $\omega_{\gamma \left( \delta _{k}\right) }\in G$
such that
\begin{equation}
   m_{\gamma \left( \delta _{k}\right) } \leq
   J_{\gamma \left( \delta
  _{k}\right) }\left( \omega_{\gamma \left( \delta _{k}\right) }\right) \leq
\frac{\delta _{k}^{2}}{2}+\frac{\gamma \left( \delta _{k}\right) }{2}\left\Vert
 \omega^{\ast }- \omega_{0}\right\Vert _{Q}^{2}.  \label{1.41}
\end{equation}

Thus, by (\ref{1.39}) and (\ref{1.41})
\begin{align}
 \frac{1}{2}\left\Vert y(\omega_{\gamma \left( \delta _{k}\right) }) -  t\right\Vert
_{B_{2}}^{2}+\frac{\gamma \left( \delta _{k}\right)}{2}   \left\Vert \omega_{\gamma \left( \delta _{k}\right) }- \omega_{0}\right\Vert _{Q}^{2} = J_{\gamma }\left( \omega_{\gamma \left( \delta _{k}\right) }  \right).
\label{tikh1}
\end{align}
From \eqref{tikh1} follows that
\begin{align}
 \frac{1}{2}\left\Vert y(\omega_{\gamma \left( \delta _{k}\right) }) -  t\right\Vert
 _{B_{2}}^{2}  \leq J_{\gamma }\left( \omega_{\gamma \left( \delta _{k}\right) }  \right), \\
\frac{\gamma \left( \delta _{k}\right)}{2}   \left\Vert \omega_{\gamma \left( \delta _{k}\right) }- \omega_{0}\right\Vert _{Q}^{2} \leq J_{\gamma }\left( \omega_{\gamma \left( \delta _{k}\right) }  \right).
\label{observe}
\end{align}
Using  \eqref{observe} and then \eqref{1.41}  one can obtain
\begin{align}
 \left\Vert \omega_{\gamma \left( \delta _{k}\right) }- \omega_{0}\right\Vert _{Q}^{2} \leq
  \frac{2}{\gamma \left( \delta _{k}\right)}  J_{\gamma }\left( \omega_{\gamma \left( \delta _{k}\right) }  \right) \leq   \frac{2}{\gamma \left( \delta _{k}\right)}\cdot
  \left[\frac{
\delta _{k}^{2}}{2}+\frac{\gamma \left( \delta _{k}\right) }{2}\left\Vert
 \omega^{\ast }- \omega_{0}\right\Vert _{Q}^{2}   \right]
\label{tikh2}
\end{align}
from what follows that
\begin{equation}
  \left\Vert \omega_{\gamma \left( \delta _{k}\right) } - \omega_{0}\right\Vert
_{Q}^{2}    \leq   \frac{\delta _{k}^{2}}{\gamma \left( \delta _{k}\right) }
+\left\Vert \omega^{\ast } - \omega_{0}\right\Vert _{Q}^{2}.  \label{1.42n}
\end{equation}
Suppose that
\begin{equation}
\lim_{k\rightarrow \infty }\gamma \left( \delta _{k}\right) =0\text{ and }%
\lim_{k\rightarrow \infty }   \frac{\delta _{k}^{2}}{\gamma \left( \delta
_{k}\right) } =0.  \label{1.43}
\end{equation}%
Then (\ref{1.42n}) implies that the sequence $\left\{ \omega_{\gamma \left( \delta
_{k}\right) }\right\} \subset G\subseteq Q$ is bounded in the norm of the
space $Q.$ Since $Q$ is compactly embedded in $B_{1},$ then there exists a
sub-sequence of the sequence $\left\{ \omega_{\gamma \left( \delta _{k}\right)
}\right\} $ which converges in the norm of the space $B_{1}.$

To ensure  \eqref{1.43}
one
can choose, for example
\begin{equation}\label{regcond}
\gamma \left( \delta _{k}\right) =C\delta _{k}^{\mu
},\mu \in \left( 0, 2\right),~~C=const. > 0, \delta \in (0,1).
\end{equation}
Other choices of $\gamma$ which satisfy  conditions \eqref{1.43} are
also possible.

In \cite{BaK, BKK} was proposed  following iterative update of the regularization parameters ${\gamma_k}$ which satisfy conditions \eqref{1.43}:
\begin{equation}\label{iterreg}
\gamma_k = \frac{\gamma_0}{(k+1)^p}, ~~ p \in (0,1],
  \end{equation}
where $\gamma_0$ can be computed as in \eqref{regcond}.

\subsection{Morozov's discrepancy principle}

\label{sec:morozov}

The principle determines the regularization parameter $\gamma=\gamma(\delta)$
  in \eqref{1.39}  such that

  \begin{equation}
\| y(\omega_{\gamma(\delta)}) - t \| = c_m \delta,
\label{tikhregfunc}
    \end{equation}
  where $c_m \geq 1$ is a constant.
Relaxed version of a  discrepancy principle is:
 \begin{equation}
c_{m,1} \delta \leq \| y(\omega_{\gamma(\delta)}) - t \|  \leq c_{m,2} \delta,
\label{relaxtikhregfunc}
    \end{equation}
  for some constants $1 \leq c_{m,1} \leq c_{m,2}$.
  The main feature of the principle is that the computed weight function
  $\omega_{\gamma(\delta)}$  can't be more accurate than the residual
  $\| y(\omega_{\gamma(\delta)}) - t \|$.

For the Tikhonov functional  
\begin{equation}
J_{\gamma}\left( \omega\right) =\frac{1}{2}\left\Vert y(\omega)- t\right\Vert
_2^{2}+  \gamma \| \omega \|_2^2 = \varphi(\omega) + \gamma \psi(\omega),
\label{balancingtikh}
\end{equation}
the value function $F(\gamma): \mathbb{R}^+ \to \mathbb{R}$
 is defined  accordingly to \cite{TA}  as
  \begin{equation}
    F(\gamma) = \inf_{\omega} J_{\gamma}(\omega).
    \label{valuefunc}
    \end{equation}

 If there exists $F_\gamma'(\gamma)$  at $\gamma >0$ then  from \eqref{balancingtikh} and
  \eqref{valuefunc}  follows that
  \begin{equation}
    F(\gamma) =\inf_{\omega} J_{\gamma }\left( \omega\right)
    = \underbrace{\varphi'(\omega)}_{\bar{\varphi}(\gamma)}
    + \gamma  \underbrace{\psi'(\omega)}_{ \bar{\psi}(\gamma)}.
\label{tikhfunc2}
  \end{equation}
  Since $ F_{\gamma}'(\gamma) = \psi'(\omega) =  \bar{\psi}(\gamma)$  then from \eqref{tikhfunc2} follows
  \begin{equation}
    \bar{\psi}(\gamma) =  F_\gamma'(\gamma),~~ \bar{\varphi}(\gamma) = F(\gamma) - \gamma F_\gamma'(\gamma).
    \end{equation}

The main idea of the principle is to compute discrepancy $\bar{\varphi}(\gamma)$  using the value function  $F(\gamma)$  and then approximate $F(\gamma)$  via  model functions.
  If  $ \bar{\psi}(\gamma) \in C(\gamma)$  then the discrepancy equation
  \eqref{tikhregfunc}
  can be rewritten  as
  \begin{equation}\label{modelfunc} 
  \bar{\varphi}(\gamma) =  F(\gamma) - \gamma F_\gamma'(\gamma) = \frac{\delta^2}{2}.
  \end{equation}
  The goal is to solve \eqref{modelfunc}  for $\gamma$. 
 Main methods for solution of \eqref{modelfunc} are the model function approach and a quasi-Newton method  presented in details in \cite{IJ}.

\subsection{Balancing principle}

\label{sec:balancing}

The balancing principle (or Lepskii, see \cite{LLP, mat})
finds   $\gamma > 0$   such that following expression is fulfilled
   \begin{equation}\label{balancing1}
     \bar{\varphi}(\gamma) = C  \gamma \bar{\psi}(\gamma),
     \end{equation}
   where $C = a_0/a_1$ is determined by the statistical a priori knowledge  from shape parameters in Gamma distributions \cite{IJ}. When $C=1$ the method is called zero crossing  method, see details in \cite{JG}.
  For   iterative update of $\gamma$   in \cite{IJ}  was  proposed  the  fixed point algorithm \ref{alg:fp}. Convergence  of this algorithm  is  also analyzed in \cite{IJ}.

\begin{algorithm}[hbt!]
  \centering
  \caption{Fixed point algorithm \label{alg:fp}}
  \begin{algorithmic}[1]
    \STATE   Start with the initial approximations
  $\gamma_{0}$    and compute the sequence of ${\gamma}_k$  in the following steps.

    \STATE  Compute the value function $F(\gamma_k) = \inf_{\omega} J_{\gamma_k}(\omega)$  for \eqref{balancingtikh}   and get $\omega_{\gamma_k}$.

   \STATE    Update the reg. parameter $\gamma :=
  \gamma_{k+1}$   as
\begin{equation*}
\begin{split}
\gamma_{k+1} &= \frac{\|\bar{\varphi}(\omega_{\gamma_k})\|_2}{\|\bar{\psi}(\omega_{\gamma_k})\|_2}
\end{split}
\end{equation*}

\STATE For the tolerance  $0 <\theta < 1$ chosen by the user, stop computing reg.parameters
  $\gamma_{k}$ if computed
  $\gamma_{k} $ are stabilized, or $|\gamma_{k} -\gamma_{k-1} | \leq \theta$.  Otherwise, set $k:=k+1$ and go to Step 2.

  \end{algorithmic}
\end{algorithm}

\begin{figure}[hbt!]
 \begin{center}
   \begin{tabular}{cc}
{\includegraphics[scale=0.4, clip=true,]{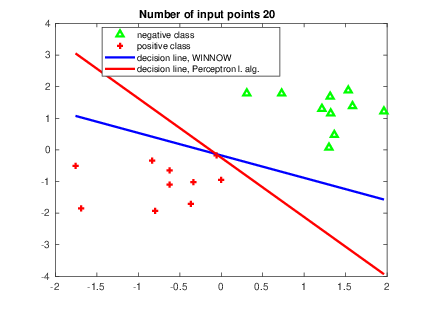}} &
{\includegraphics[scale=0.4,  clip=true,]{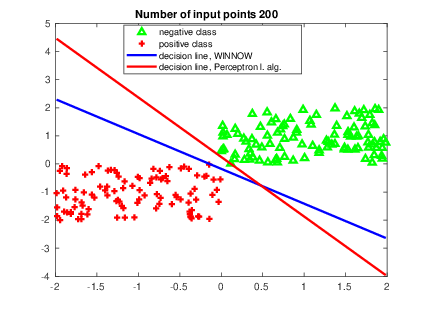}}  
\end{tabular}
 \end{center}
 \caption{Comparison of two classification algorithms for separation
   of two classes: Perceptron learning algorithm (red line) and WINNOW
   (blue line).}
  \label{fig:6}
\end{figure}

\section{Numerical   results}

\label{sec:numex}

In this section are presented several examples which show performance
and effectiveness of least squares, perceptron and WINNOW algorithms
for classification. We note that all classification algorithms
considered here doesn't include regularization.

\subsection{Test 1}

In this test the goal is to compute decision boundaries for two linearly
separated classes using least squares classification. Points in these classes are generated randomly
 by
the   linear function $y = 1.2 - 0.5x$ 
 on different input intervals for $x$. Then the random noise $\delta$  is added to the data $y(x)$  as
\begin{equation}\label{noise}
    y_\delta(x) = y(x)(1 + \delta \alpha),
\end{equation}
where $\alpha \in (−1, 1)$ is randomly distributed number and $\delta
\in [0, 1]$ is the noise level.
Then obtained points are classified such that the target function  for classification is defined as
\begin{equation}\label{test1}
    t_i = \left \{
    \begin{array}{ll}
     1, &   \textrm{$y_\delta(x) - 1.2 + 0.5x > 0$}, \\
     0,  &  \textrm{otherwise}. 
    \end{array}
    \right.
    \end{equation}

Figures \ref{fig:1} present
classification performed via least squares minimization $\min_\omega \| A\omega
-y \|_2^2$ for the linear model function \eqref{plmodel0} with target
values $t$ given by \eqref{test1} and elements of the design matrix
$A$ given by \eqref{plmodel2}. 
Using these figures we observe that the least squares can be used very
successfully for classification when it is known that classes are linearly
separated.

\subsection{Test 2}

Here we present examples of performance of the perceptron learning
algorithm for classification of two linearly separated classes. Data
for analysis in these examples are generated similarly as in the Test
1.  Figures \ref{fig:4} present classification of two classes
in the perceptron algorithm with three
weights. Figures \ref{fig:5} show classification of two classes using
the second order polynomial function \eqref{secondorder} in the
perceptron algorithm with six weights. We note that the red and blue
lines presented in Figures \ref{fig:5} are classification boundaries computed via
\eqref{secondorder6}. Figures \ref{fig:6} present comparison of linear
perceptron and WINNOW algorithms for separation of two classes.
Again, all these figures show that perceptron and WINNOW algorithms
can be successfully  used for separation of linearly separated classes.

\begin{figure}[hbt!]
 \begin{center}
   \begin{tabular}{c}
      {\includegraphics[scale=0.35, clip=true,]{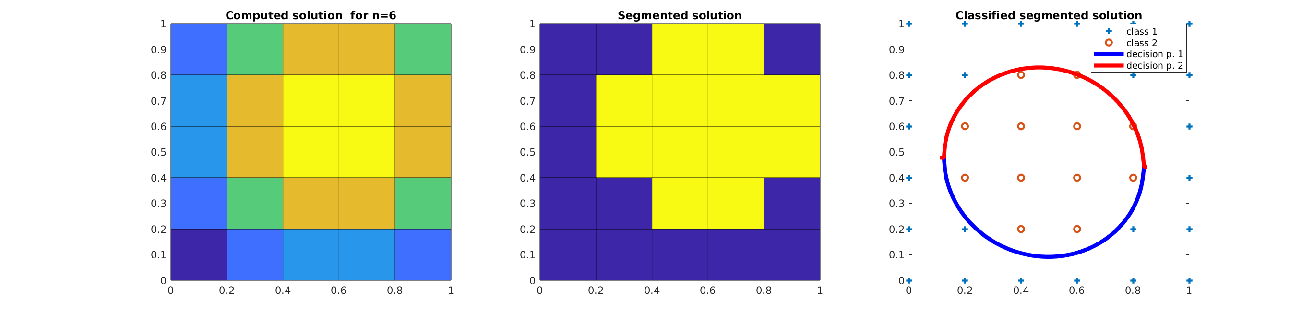}}  \\
     {\includegraphics[scale=0.35, clip=true,]{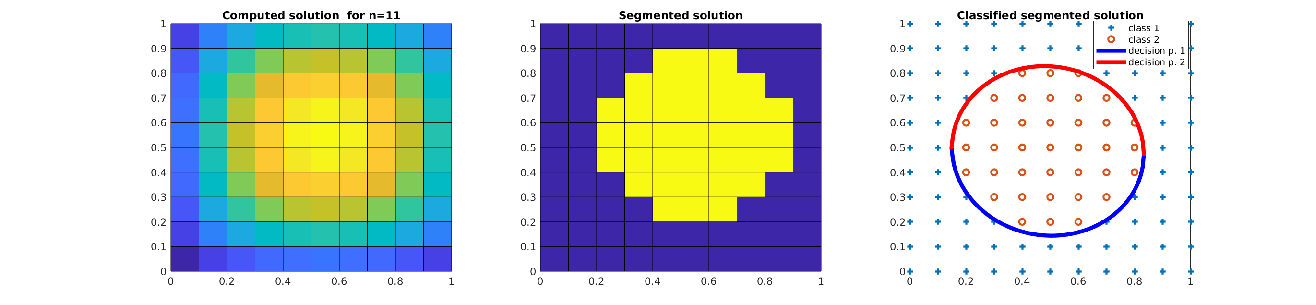}}  \\
      {\includegraphics[scale=0.35, clip=true,]{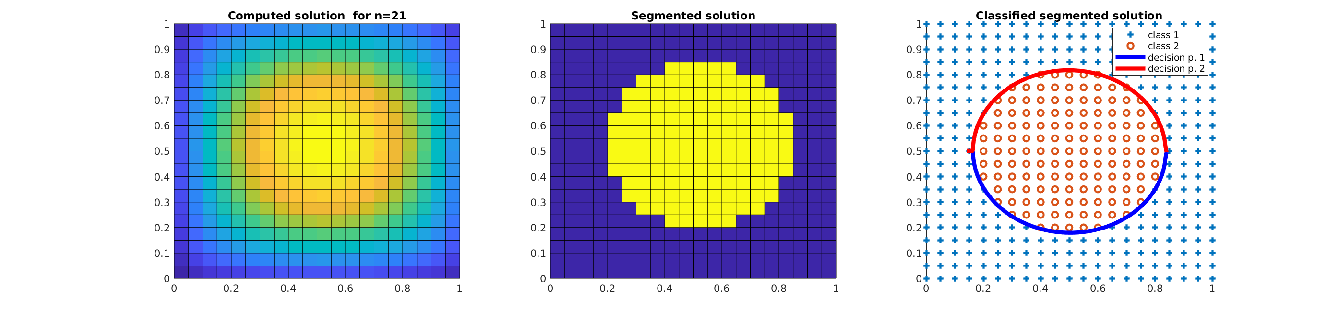}}   
\end{tabular}
 \end{center}
 \caption{Classification of the computed solution for Poisson's equation on the unit square (see example 8.1.3 of \cite{BKK})  for different number of mesh points $n^2$ using perceptron learning algorithm. }
  \label{fig:7}
\end{figure}

\subsection{Test 3}

This test shows performance of using the second order polynomial
function \eqref{secondorder} in the perceptron algorithm for  classification of the
segmented solution.

The   classification  problem in this example is as
follows:

\begin{itemize}

  \item  Given the computed solution of the Poisson's equation
$\triangle u = f$ with homogeneous boundary conditions $u=0$ on the
    unit square, classify the  discrete solution $u_h$  into two classes
such that the  target function for classification is  defined as (see the top right figure  of  Figure \ref{fig:7}):
\begin{equation}\label{test3_1}
    t_i = \left \{
    \begin{array}{lll}
     1, &    \textrm{${u_h}_i  > 4$} & \textrm{(yellow points)}, \\
     0,  &  \textrm{otherwise} &  \textrm{(blue  points)}. \\
    \end{array}
    \right.
    \end{equation}

\item   Use the second order polynomial function \eqref{secondorder} in the
perceptron algorithm to compute decision boundaries.

  \end{itemize}

For details about setup
of the computed solution for Poisson's equation we refer to the
example 8.1.3 of \cite{BKK}.  Figure \ref{fig:7} presents results of
using the second order polynomial for the classification of the computed segmented solution. The
computed solution $u_h$ of the Poisson's equation on different meshes is presented on the
left figures of Figure \ref{fig:7}. The middle figures show segmentation of the
obtained solution  satisfying \eqref{test3_1}.
The  right figures of Figure \ref{fig:7} show results of applying
the second order polynomial function \eqref{secondorder} in the
perceptron algorithm for classification of the computed solution $u_h$ with
target function \eqref{test3_1}. We observe, that
computed decision points  correctly separates two classes
 even if these classes are not separable.

\subsection{Test 4}

\begin{figure}[hbt!]
\begin{center}
\begin{tabular}{cc}
{\includegraphics[width = 0.5 \textwidth]{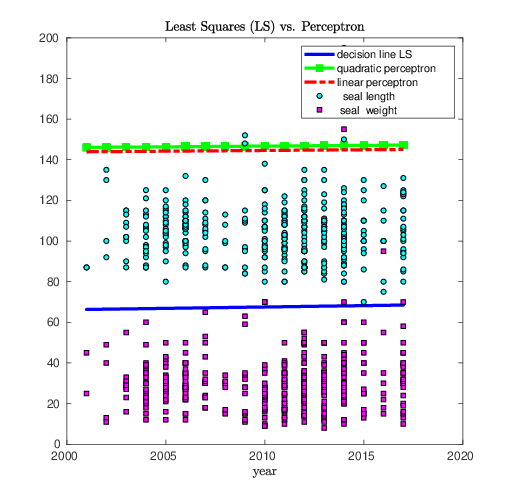}} &
{\includegraphics[width = 0.5 \textwidth]{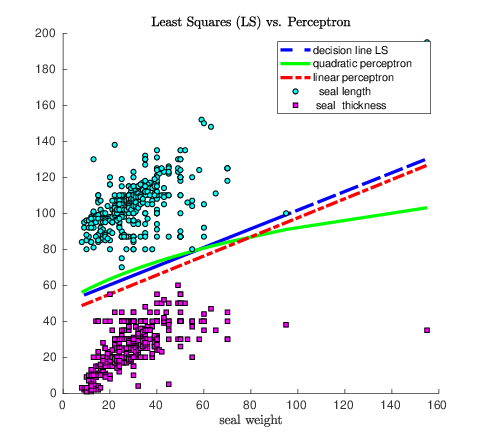}}
\\
a) & b)
\end{tabular}
\end{center}
\caption{Comparison of least squares (LS) and Perceptron learning algorithm for separation of two classes using Grey Seal database \cite{waves}.}
  \label{fig:8}
\end{figure}

In this test we show performance of linear least squares together with
linear and quadratic perceptron algorithms for classification of
experimental data sets: database of grey seals \cite{waves} and Iris
data set \cite{iris}.  Figure \ref{fig:8}-a) shows classification of
seal length and seal weight depending on the year.  Figure
\ref{fig:8}-b) shows classification of seal length and seal thickness
depending on the seal weight. We observe that classes on Figure
\ref{fig:8}-a) are not linearly separable and the best algorithm which
separates both classes well, is the least squares. In this example,
the linear and quadratic perceptron have not classified data correctly
and actually, these algorithms have not converged and stopped when the
maximal number of iterations ($10^8$) was reached.
As soon as classes become  linearly separated, all algorithms show good performance and  computes almost the same separation lines, see Figure \ref{fig:8}-b).

\begin{figure}[hbt!]
\begin{center}
\begin{tabular}{c}
{\includegraphics[width = 0.5  \textwidth]{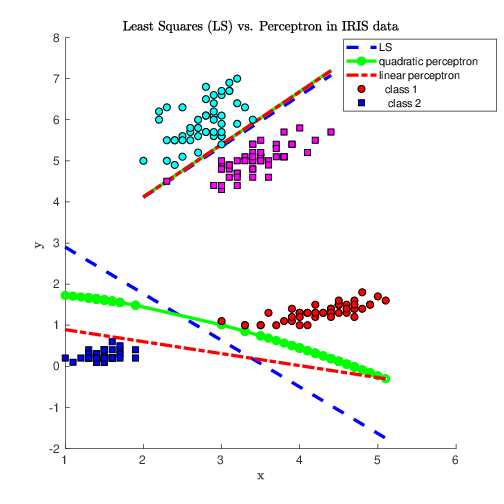}}
\end{tabular}
\end{center}
\caption{Comparison of least squares (LS) and Perceptron learning algorithm on Iris dataset \cite{iris}.}
  \label{fig:9}
\end{figure}

The same conclusion is obtained from separation of Iris data set
\cite{iris}.  Figures \ref{fig:9} show decision lines computed by
least squares, linear and quadratic perceptron algorithms. Since all
classes of Iris data set are linearly separable, all classification algorithms
separate data correctly.

\section{Conclusions}

We have presented regularized and non-regularized perceptron learning
and least squares algorithms for classification problems as well as
discussed main a-priori and a-posteriori Tikhonov's regularization
  rules for choosing the regularization parameter.
 The Fr\'{e}chet derivatives for
 least squares  and perceptron algorithms are  also rigorously
  derived.

The future work can be related to computation of miss-classification
rates which can be done similarly with works \cite{hand, thomas}, as
well as to study of classification problem using regularized linear
regression, perceptron learning and WINNOW algorithms. Other
classification algorithms such that regularized SVM and kernel methods
can be also analyzed. Testing of all algorithms on different benchmarks
as well as extension to multiclass case should be also investigated.

\section*{Acknowledgments}

 The research 
is supported by the Swedish Research Council grant VR 2018-03661.

\end{document}